\documentclass[preprint,12pt]{elsarticle}

\usepackage{amsmath,amssymb,amscd}
\usepackage{amsthm}
\usepackage{subfigure}
\usepackage[nodots]{numcompress}
\usepackage{verbatim}
\usepackage{fancybox}
\usepackage{subfigure}
\usepackage[usenames,dvipsnames]{color}
\usepackage{sidecap}
\usepackage{caption}
\usepackage{ifthen}
\usepackage{textcomp}
\usepackage[all]{xy}

\newcommand{\Q}{{\mathbb Q}}

\newcommand{\dem}{{\em Proof: \;}}
\newcommand{\fdem}{\hfill $\square$}

\theoremstyle{plain}
\newtheorem{teo}{Theorem}[section]
\newtheorem{lema}[teo]{Lemma}
\newtheorem{cor}[teo]{Corollary}
\newtheorem{prop}[teo]{Proposition}

\theoremstyle{definition}
\newtheorem{defi}{Definition}[section]

\theoremstyle{remark}

\begin{document}

\begin{frontmatter}

\title{Factorization of network reliability with perfect nodes II: Connectivity matrix}

\author[IMATE]{Juan Manuel Burgos}
\address[IMATE]{Instituto de Matem\'aticas, Universidad Nacional Aut\'onoma de M\'exico, Unidad Cuernavaca. Av. Universidad s/n, Col. Lomas de
Chamilpa. Cuernavaca, Morelos M\'exico, 62209.\\ \texttt{\scriptsize{Email: burgos@matcuer.unam.mx}}}


\begin{abstract}
We prove the determinant connectivity matrix formula. Mathematically, the proof introduces novel techniques based on an algebraic approach and connectivity properties. Although this is the second part of a previous paper and has its original motivation there, the paper is self contained and the result is  interesting in itself.
\end{abstract}

\begin{keyword}
Network Reliability \sep Graph Theory \sep Factorization
\end{keyword}

\end{frontmatter}


\section{Introduction}

Denote by $Part_{n}$ the set of partitions of $\{1,2,\ldots n\}$. We will call the partition $\{ \{ 1,2,\ldots n \} \}$ the trivial partition. The set $Part_{n}$ has a monoid structure with unit $\{ \{1 \} , \{ 2 \} , \ldots \{ n \} \}$ under the following product: given the partitions $\mathcal{A}$ and $\mathcal{B}$, the product $\mathcal{A}\cdot \mathcal{B}$ is the finer partition that is coarser to $\mathcal{A}$ and coarser to $\mathcal{B}$. Observe that the product of any partition with the trivial one is trivial.

\begin{figure}
\begin{center}
  \includegraphics[width=0.6\textwidth]{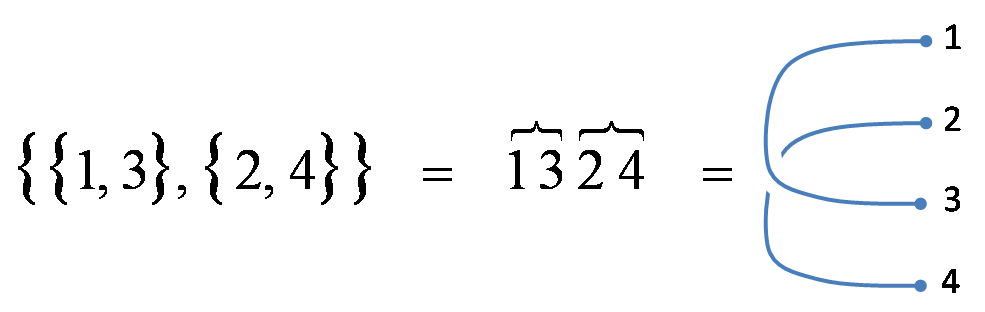}\\
  \end{center}
  \caption{}\label{EjEstadoConect}
\end{figure}

\begin{figure}
\begin{center}
  \includegraphics[width=0.5\textwidth]{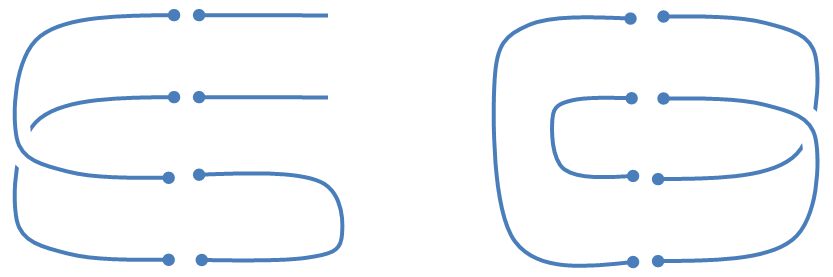}\\
  \end{center}
  \caption{Connected diagrams}\label{EjEstadoConexo}
\end{figure}

\begin{figure}
\begin{center}
  \includegraphics[width=0.5\textwidth]{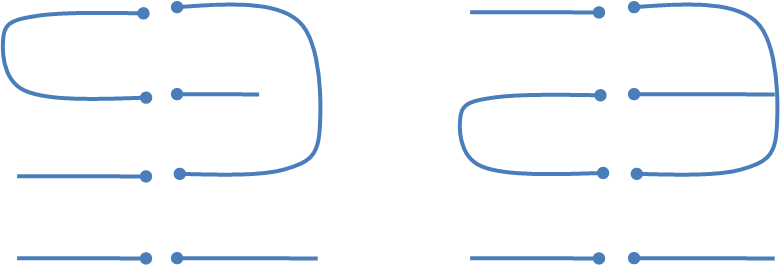}\\
  \end{center}
  \caption{Non connected diagrams}\label{EjEstadoNoConexo}
\end{figure}

Figure \ref{EjEstadoConect} shows some useful notational and diagrammatical  ways to represent a partition. Observe that the product of two partitions is trivial if and only if the resulting diagram after joining their respective vertices is connected (in the usual topological sense). Figures \ref{EjEstadoConexo} and \ref{EjEstadoNoConexo} show examples of connected and non connected diagrams.

\begin{defi}
Considering an ordering in $Part_{n}$, we define the matrix $A=(a_{ij})$ by $a_{ij}=1$ if $\mathcal{A}_{i} \cdot \mathcal{A}_{j}$ is trivial and $a_{ij}=0$ if it is not. The matrix $A$ will be called the \emph{connectivity matrix}.
\end{defi}



This paper is about the connectivity matrix and the proof that it is invertible. Moreover, we prove the main theorem of this paper:

\begin{teo}
The determinant of the connectivity matrix $A$ is: $$det(A)= \pm \prod_{\mathcal{A}\in Part_{n}}(m_{\mathcal{A}}-1)!$$ where $m_{\mathcal{A}}$ is the number of classes in the partition $\mathcal{A}$. In particular, $A$ is invertible.
\end{teo}

We observe that the left hand side of the equality is in some sense topological (it is related to connectedness) while the right hand side is purely combinatorial.

\section{Algebra of partitions}

We will identify a partition with its symmetry: i.e. there is a one to one correspondence between partitions of $\{1,2,\ldots n\}$ and the subgroups generated by transpositions of the permutation group $S_{n}$ (see [BM], [Ro]).

For this purpose, we define the orbit map $\mathcal{O}:Subgroups(S_{n})\rightarrow Part_{n}$ such that $\mathcal{O}(\mathcal{G})$ is the space of orbits of the action of $\mathcal{G}$ induced by $S_{n}$ on $\{1, 2,\ldots n\}$; i.e.  $$\mathcal{O}(\mathcal{G})= \{\ \mathcal{O}(i)\ /\ i=1,2,\ldots n\ \}= \{1, 2,\ldots n\}/\mathcal{G}$$ where $\mathcal{O}(i)$ is the set of elements $j$ for which there is a group element $g\in \mathcal{G}$ such that $g(i)=j$; i.e $\mathcal{O}(i)= \mathcal{G}\cdot i$.

On the other side, given a partition $\mathcal{A}= \{A_{1},A_{2},\ldots A_{l}\}$ we define the group $\mathcal{G}(\mathcal{A})$ generated by all transpositions $(i,j)$ for which there is some $A_{k}$ such that $i,j\in A_{k}$ and $i \neq j$. By definition, is clear that $$\mathcal{O}(\mathcal{G}(\mathcal{A}))= \mathcal{A}$$ showing in particular that $\mathcal{O}$ is surjective. As an example, consider the partition $\mathcal{A}=\{\{1,2\},\{3,4,5\},\{6\}\}$. Then, $\mathcal{G}(\mathcal{A})= \langle (1,2),(3,4),(3,5),(4,5)\rangle$. However, the orbit space of the group $\langle (1,2),(3,4,5)\rangle$ is also the partition $\mathcal{A}$ and this group is a proper subgroup of $\mathcal{G}(\mathcal{A})$. In particular, the orbit map is not injective. Nevertheless, restricted the subgroups generated by transpositions, the orbit map become one to one. This is the content of the following lemma:

\begin{lema}
Consider a partition $\mathcal{A}\in Part_{n}$ and a subgroup $\mathcal{G}< S_{n}$ such that $\mathcal{O}(\mathcal{G})= \mathcal{A}$. Then $\mathcal{G}$ is a subgroup of $\mathcal{G}(\mathcal{A})$. Moreover, if $\mathcal{G}$ is generated by transpositions, then $\mathcal{G}= \mathcal{G}(\mathcal{A})$.
\end{lema}
\dem
Consider $g\in \mathcal{G} < S_{n}$ and the partition $\mathcal{A}= \{A_{1},A_{2},\ldots A_{l}\}$. $g$ can be written as a product of disjoint cycles and because $\mathcal{O}(\mathcal{G})= \mathcal{O}(\mathcal{G}(\mathcal{A}))$, members of the same disjoint cycle belong to the same partition class; i.e. if $$g= (i_{1}^{1}\ldots i_{j_{1}}^{1})\ldots (i_{1}^{k}\ldots i_{j_{k}}^{k})$$ then for each $l$ there is a partition class $A_{m}$ such that $$i_{1}^{l},\ldots i_{j_{l}}^{l}\in A_{m}$$
By definition of the group $\mathcal{G}(\mathcal{A})$, $$(i_{1}^{l},i_{2}^{l}),\ldots (i_{1}^{l},i_{j_{l}}^{l})\in \mathcal{G}(\mathcal{A})$$ and because every disjoint cycle can be written as $$(i_{1}^{l}\ldots i_{j_{l}}^{l})=  (i_{1}^{l},i_{j_{l}}^{l})\ldots (i_{1}^{l},i_{3}^{l})(i_{1}^{l},i_{2}^{l})$$
we conclude that $g$ belongs to $\mathcal{G}(\mathcal{A})$. This way $\mathcal{G}<\mathcal{G}(\mathcal{A})$.

Now assume that $\mathcal{G}$ is generated by transpositions. Because $(i,k)=(j,k)(i,j)(j,k)$ we have that $(i,j),(j,k)\in \mathcal{G}$ implies $(i,k)\in \mathcal{G}$ and the following equivalence relation: $i\sim j$ are equivalent if $i=j$ or $(i,j)\in \mathcal{G}$. The quotient space of this equivalence relation defines a partition $$\mathcal{A}'= \{1,2,\ldots n\}/\sim$$
Because $\mathcal{G}$ is generated by transpositions, we have $\mathcal{A}'= \mathcal{O}(\mathcal{G})=\mathcal{A}$ and by definition of $\mathcal{G}(\mathcal{A})$ and the equivalence relation, we have that $\mathcal{G}(\mathcal{A})< \mathcal{G}$. Because $\mathcal{G}(\mathcal{A})> \mathcal{G}$, we conclude that $\mathcal{G}(\mathcal{A})= \mathcal{G}$.
\fdem

\begin{cor}
Restricted to the subgroups generated by transpositions, the orbit map is one to one.
\end{cor}

Moreover, under the product $$\mathcal{G}_{1}\cdot\mathcal{G}_{2}= \langle\mathcal{G}_{1}\cup\mathcal{G}_{2}\rangle$$ the orbit map restricted to the subgroups generated by transpositions is a monoid isomorphism: $\mathcal{O}(\mathcal{G}_{1}\cdot\mathcal{G}_{2})= \mathcal{O}(\mathcal{G}_{1})\cdot \mathcal{O}(\mathcal{G}_{2})$ and $\mathcal{O}( \{id\}) = \{\{1\}, \{2\},\ldots \{n\}\}$.

Because the product of subgroups generated by transpositions is generated by transpositions, this product is well defined. Observe that $\mathcal{O}(\mathcal{G}_{1})$ is finer than $\mathcal{O}(\mathcal{G}_{2})$ as partitions if and only if $\mathcal{G}_{1}\subset \mathcal{G}_{2}$ as subgroups of $S_{n}$.

From now on we will identify partitions of $\{1,2,\ldots n\}$ with subgroups generated by transpositions of the permutation group $S_{n}$ and call these elements connectivity states. Under this identification we rewrite the definition of the connectivity matrix as follows:

\begin{defi}\label{ConnMatrix}
Relative to an ordering in $Part_{n}$ the connectivity matrix $A=(a_{ij})$ is the matrix given by $a_{ij}=1$ if $\mathcal{A}_{i}\cdot \mathcal{A}_{j}= S_{n}$ and $a_{ij}=0$ otherwise.
\end{defi}

Observe that every partition $\mathcal{A}$ is idempotent ($\mathcal{A}^{2}= \mathcal{A}$) and if $\mathcal{A}<\mathcal{B}$ then $\mathcal{A}\cdot\mathcal{B}= \mathcal{B}$.

\begin{defi}
The permutation group $S_{n}$ acts by conjugation in $Part_{n}$: $$\sigma(\mathcal{A})= \sigma \mathcal{A} \sigma^{-1}$$ where $\mathcal{A}\in Part_{n}$ and $\sigma \in S_{n}$.
\end{defi}

In particular we have conjugation classes in $Part_{n}$; i.e. the orbits in $Part_{n}$ under the action of $S_{n}$: $\mathcal{A}\sim \mathcal{B}$ if there is a permutation $\sigma$ such that $\mathcal{A}=\sigma(\mathcal{B})=\sigma \mathcal{B}\sigma^{-1}$.

\begin{lema}
Consider a permutation $\sigma$ in $S_{n}$. Then $\sigma:Part_{n}\rightarrow Part_{n}$ is a monoid morphism; i.e. $\sigma(\{id\})=\{id\}$ and $$\sigma(\mathcal{A}\cdot\mathcal{B})=\sigma(\mathcal{A})\cdot\sigma(\mathcal{B})$$
\end{lema}
\dem
It is clear that $\sigma(\{id\})=\{id\}$. Let's show that $\sigma$ is a monoid morphism. Observe that $\sigma(\mathcal{A}\cdot\mathcal{B})$ contains $\sigma(\mathcal{A})$ and $\sigma(\mathcal{B})$ so, by definition, $$\sigma(\mathcal{A}\cdot\mathcal{B})\supset\sigma(\mathcal{A})\cdot\sigma(\mathcal{B})$$ See also that, $\sigma^{-1}(\sigma(\mathcal{A})\cdot\sigma(\mathcal{B}))$ contains $\mathcal{A}$ and $\mathcal{B}$ so, by definition, $$\mathcal{A}\cdot\mathcal{B}\subset\sigma^{-1}(\sigma(\mathcal{A})\cdot\sigma(\mathcal{B}))$$ Acting by $\sigma$ we get $$\sigma(\mathcal{A}\cdot\mathcal{B})\subset\sigma(\mathcal{A})\cdot\sigma(\mathcal{B})$$
\fdem

Because $Part_{n}$ is a commutative monoid, the $\Q$-vector space $V$ generated by $Part_{n}$ introduced before is actually an associative and commutative $\Q$-algebra with unit. We call this algebra the algebra of partitions. The monoid morphism $\sigma$ extends linearly to a unital algebra endomorphism where $\sigma$ is a permutation.

\section{The determinant formula}

\subsection{Introduction}

This section is devoted to the calculation of the connectivity matrix determinant and to show in particular that this matrix is invertible. In general terms, the proof goes in three steps:

\begin{enumerate}
  \item Gauss elimination: Through elementary row operations, for each $n$ we develop a Gauss elimination method such that the connectivity matrix becomes lower triangular.
  \item Connectivity numbers and reliability polynomial: Identify the resulting diagonal elements, called connectivity numbers, with the coefficient of the higher order term of certain reliability polynomials.
  \item Connectivity numbers calculation: Under these identification, calculate the diagonal elements.
\end{enumerate}

\subsection{Step one: Gauss elimination}

To develop a Gauss elimination on the connectivity matrix, we must develop a systematic method to detect whether the entries of the matrix are zero or one. After a careful analysis, one notice that under particular orderings of the $Part_{n}$ basis (these will be called coherent orderings later), some blocks related to the symmetry of the partitions appear. This suggests that it is unavoidable the consideration of the partition's symmetry in order to develop a Gauss elimination on the connectivity matrix. The partition's symmetry is encoded in the partitions algebra introduced in the previous section. Non trivial identities of this algebra will give the Gauss elimination we are looking for.

Consider the linear operator $\pi:V\rightarrow V$ such that $\pi(S_{n})=S_{n}$ and for every connectivity state $\mathcal{A}$ distinct from $S_{n}$, $$\pi(\mathcal{A})=\prod_{\mathcal{B}\ /\ \mathcal{B}\nless \mathcal{A}}(\mathcal{A}-\mathcal{A}\cdot\mathcal{B})$$

\begin{lema}\label{lemaNueve}
Consider a connectivity state $\mathcal{A}$. The vector $\pi(\mathcal{A})$ satisfies the following properties:

\begin{enumerate}
  \item $\mathcal{B}\cdot\pi(\mathcal{A})=0$   $\forall\mathcal{B}\ /\ \mathcal{B}\nless \mathcal{A}$
  \item $\mathcal{C}\cdot\pi(\mathcal{A})=\pi(\mathcal{A})$   $\forall\mathcal{C}\ /\ \mathcal{C}\leq \mathcal{A}$
  \item $\pi(\sigma(\mathcal{A}))=\sigma(\pi(\mathcal{A}))$   $\forall\sigma\in S_{n}$
\end{enumerate}

In particular, $\mathcal{A}\cdot\pi(\mathcal{A})=\pi(\mathcal{A})$ and $\mathcal{B}\cdot\pi(\mathcal{A})=0$ for every connectivity state $\mathcal{B}$ distinct and conjugated to $\mathcal{A}$.
\end{lema}
\dem
\begin{enumerate}
  \item The algebra is commutative and $Part_{n}$ is a basis of idempotents; i.e. $\mathcal{A}^{2}=\mathcal{A}$ for every connectivity state $\mathcal{A}$. In particular, there is a factor of $\mathcal{B}\cdot\pi(\mathcal{A})$ that is zero:
  $$\mathcal{B}\cdot(\mathcal{A}-\mathcal{A}\cdot\mathcal{B})= \mathcal{B}\cdot\mathcal{A}-\mathcal{B}^{2}\cdot\mathcal{A}=0$$
  \item Because of the fact that $\mathcal{C}\cdot\mathcal{A}=\mathcal{A}$ if $\mathcal{C}\leq \mathcal{A}$, every factor of $\pi(\mathcal{A})$ remains the same after multiplying by $\mathcal{C}$:
  $$\mathcal{C}\cdot(\mathcal{A}-\mathcal{A}\cdot\mathcal{B})= \mathcal{C}\cdot\mathcal{A}-\mathcal{C}\cdot\mathcal{A}\cdot\mathcal{B}= \mathcal{A}-\mathcal{A}\cdot\mathcal{B}$$
  \item \begin{eqnarray*}
          \sigma(\pi(\mathcal{A})) &=& \sigma\left(\prod_{\mathcal{B}\ /\ \mathcal{B}\nless \mathcal{A}}(\mathcal{A}-\mathcal{A}\cdot\mathcal{B})\right) \\
           &=& \prod_{\mathcal{B}\ /\ \mathcal{B}\nless \mathcal{A}}(\sigma(\mathcal{A})-\sigma(\mathcal{A})\cdot\sigma(\mathcal{B})) \\
           &=& \prod_{\sigma(\mathcal{B})\ /\ \sigma(\mathcal{B})\nless \sigma(\mathcal{A})}(\sigma(\mathcal{A})-\sigma(\mathcal{A})\cdot\sigma(\mathcal{B})) \\
           &=& \pi(\sigma(\mathcal{A}))
        \end{eqnarray*} where we used in the last identity that $\sigma$ is invertible.
\end{enumerate}
\fdem

\begin{defi}
For each connectivity state $\mathcal{A}$ we define its connectivity number as the coefficient of $S_{n}$ in the expansion of $\pi(\mathcal{A})$ relative to the basis $Part_{n}$; i.e. $$\pi(\mathcal{A})= \mathcal{A} + \ldots + \alpha_{\mathcal{A}}S_{n}$$
\end{defi}

Observe that, by the third item of Lemma \ref{lemaNueve}, the connectivity number is invariant under conjugation: $$\alpha_{\mathcal{A}}=\alpha_{\sigma(\mathcal{A})}$$ for every permutation $\sigma$.

\begin{defi}
An ordering of the basis $Part_{n}$ will be called coherent if it satisfies $$\mathcal{A}_{i} < \mathcal{A}_{j}\ \Rightarrow\ i<j$$
\end{defi}

We argue that a coherent ordering always exist in the following way: Consider the Hasse diagram (partial ordering diagram) of connectivity states in $Part_{n}$. Because conjugated states necessary belong to the same level of the Hasse diagram, we can order $Part_{n}$ in the following way: We order some conjugation class $\mathcal{O}_{i}$ of the first level, then we order some other conjugation class $\mathcal{O}_{j}$ of the same level and we continue the process until we have order all the conjugation classes of the first level. After that, we order the second level in the same way as we did in the first and so on until we have order all the partitions. The previous argument is formalized in the next lemma:

\begin{lema}
The partial ordering on $Part_{n}$ induces a partial ordering on the conjugation classes.
\end{lema}
\dem
Define the following partial order on the conjugation classes: $\mathcal{O}_{i}< \mathcal{O}_{j}$ if there is $\mathcal{A}_{i}\in\mathcal{O}_{i}$ and $\mathcal{A}_{j}\in \mathcal{O}_{j}$ such that $\mathcal{A}_{i}< \mathcal{A}_{j}$. This partial order relation is well defined because of the following fact: Suppose there are $\mathcal{A}_{i},\ \mathcal{A}'_{i}\in\mathcal{O}_{i}$ and $\mathcal{A}_{j},\ \mathcal{A}'_{j}\in\mathcal{O}_{j}$ such that $\mathcal{A}_{i}< \mathcal{A}_{j}$ and $\mathcal{A}'_{i}> \mathcal{A}'_{j}$. There are permutations $\sigma,\ \eta\in S_{n}$ such that $\mathcal{A}'_{i}= \sigma(\mathcal{A}_{i})$ and $\mathcal{A}'_{j}= \eta(\mathcal{A}_{j})$ and we have $$\mathcal{A}_{i}<\mathcal{A}_{j}= \eta^{-1}(\mathcal{A}'_{j})< \eta^{-1}(\mathcal{A}'_{i})= \eta^{-1}\sigma(\mathcal{A}_{i})$$ We conclude that $\eta^{-1}\sigma= id$ and $\mathcal{A}_{i}= \mathcal{A}_{j}$ which is absurd.
\fdem

As an example consider the $n=4$ case. The conjugation classes are:

\begin{eqnarray*}
  \mathcal{O}_{1} &=& \{\ 1234\ \} \\
  \mathcal{O}_{2} &=& \{\ 12 \overbrace{34}, 13\overbrace{24}, 23\overbrace{14}, 1\overbrace{23}4, \overbrace{13}24, \overbrace{12}34\ \} \\
  \mathcal{O}_{3} &=& \{\ \overbrace{14}\overbrace{23}, \overbrace{13}\overbrace{24}, \overbrace{12}\overbrace{34}\ \} \\
  \mathcal{O}_{4} &=& \{\ 1\overbrace{234}, 2\overbrace{134}, 3\overbrace{124},\overbrace{123}4\ \} \\
  \mathcal{O}_{5} &=& \{\ \overbrace{1234}\ \}
\end{eqnarray*}

\begin{figure}
\begin{center}
  \includegraphics[width=0.25\textwidth]{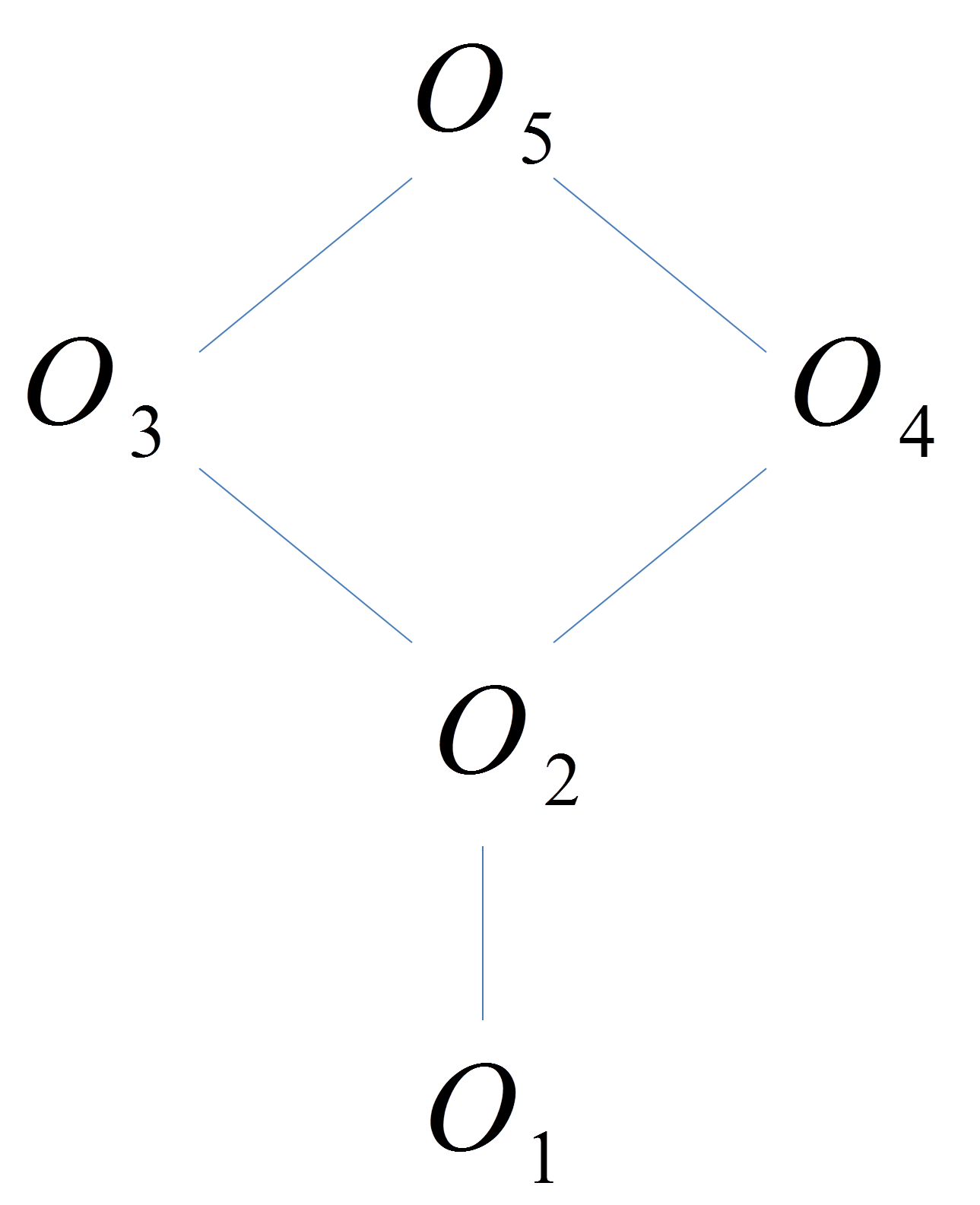}\\
  \end{center}
  \caption{Hasse diagram of conjugation classes $Part_{4}$}\label{Hasse4}
\end{figure}

Figure \ref{Hasse4} shows the Hasse diagram of the induced partial ordering on the conjugation classes of $Part_{4}$ described in the previous lemma. Is clear then that the following is a coherent order of $Part_{4}$ (we use the standard linear algebra abuse of set notation concerning ordered basis):
$$Part_{4}= \mathcal{O}_{1}\cup \mathcal{O}_{2}\cup \mathcal{O}_{3}\cup \mathcal{O}_{4}\cup \mathcal{O}_{5}$$

Observe that, once we have order the basis $Part_{n}$ in a coherent way, the operator $\pi$ is the Gauss elimination we were looking for. In fact, $\pi(\mathcal{A})= \mathcal{A} + \ldots + \alpha_{\mathcal{A}}S_{n}$ are the elementary operations of adding multiples of rows $\mathcal{B}>\mathcal{A}$ to the row $\mathcal{A}$. Because of definition \ref{ConnMatrix}: $a_{ij}=1$ if and only if $\mathcal{A}_{i}\cdot \mathcal{A}_{j}= S_{n}$ and $a_{ij}=0$ otherwise, we conclude that after the elementary operations described by $\pi(\mathcal{A}_{i})$ on the $i$-th row, the resulting entries $a'_{ij}$ equals the coefficient of $S_{n}$ in the expansion of $\pi(\mathcal{A}_{i})\cdot \mathcal{A}_{j}$. By Lemma \ref{lemaNueve} we have that $a'_{ii}= \alpha_{\mathcal{A}_{i}}$ and $a'_{ij}=0$ if $j>i$. We conclude that performing the elementary row operations $\pi(\mathcal{A}_{i})$ on the row $i$-th row on every row of the matrix, the resulting matrix is lower triangular.

\begin{figure}
\begin{center}
  \includegraphics[width=0.35\textwidth]{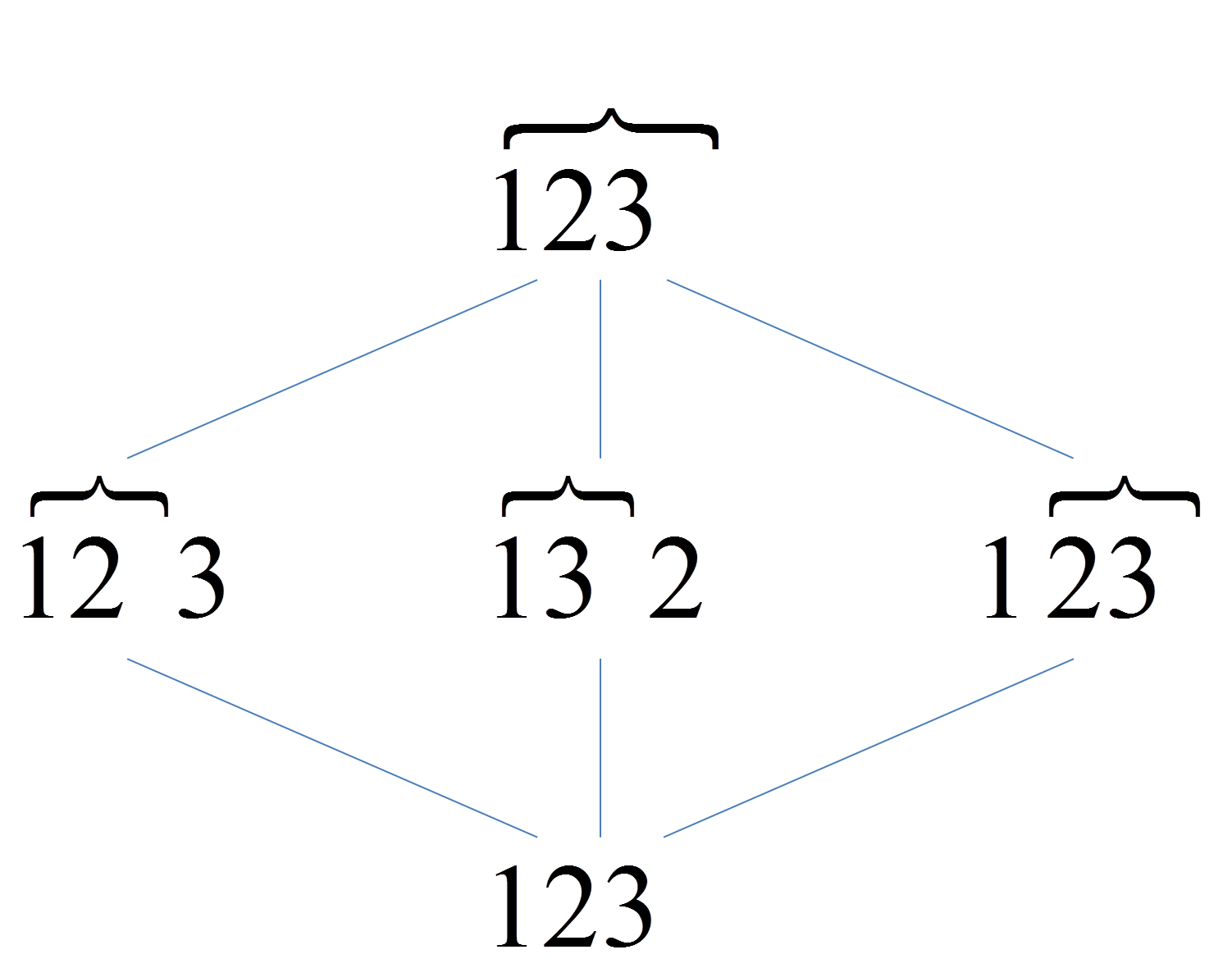}\\
  \end{center}
  \caption{Hasse diagram of $Part_{3}$}\label{Hasse3}
\end{figure}

As an example consider the $n=3$ case. The Hasse diagram of $Part_{3}$ is shown in Figure \ref{Hasse3} and is clear that a coherent order on the basis $Part_{3}$ is $$Part_{3}= \{123, 1\overbrace{23}, \overbrace{13}2, \overbrace{12}3, \overbrace{123}\}$$
and the resulting connectivity matrix is:
$$A=\left(
                                                                \begin{array}{ccccc}
                                                                  0 & 0 & 0 & 0 & 1 \\
                                                                  0 & 0 & 1 & 1 & 1 \\
                                                                  0 & 1 & 0 & 1 & 1 \\
                                                                  0 & 1 & 1 & 0 & 1 \\
                                                                  1 & 1 & 1 & 1 & 1 \\
                                                                \end{array}
                                                              \right)$$

The linear operator $\pi:V\rightarrow V$ reads as follows:

\begin{eqnarray*}
  \pi(\ \  123\ \ \ ) &=& 123\ -\ \overbrace{12}3\ -\ \overbrace{13}2\ -\ 1\overbrace{23}\ +\ 2.\ \overbrace{123}\  \\
  \pi(\ \overbrace{12}3\ ) &=& \overbrace{12}3\ -\ \overbrace{123}\  \\
  \pi(\ \overbrace{13}2\ ) &=& \overbrace{13}2\ -\ \overbrace{123}\  \\
  \pi(\ 1\overbrace{23}\ ) &=& 1\overbrace{23}\ -\ \overbrace{123}\  \\
  \pi(\ \ \overbrace{123}\ \ ) &=& \overbrace{123}
\end{eqnarray*} and the matrix associated to $\pi$ on the basis $Part_{3}$ is $$B=\left(\begin{array}{ccccc} 1 &  &  &  &  \\ -1 & 1 &  &  &  \\ -1 & 0 & 1 &  &  \\ -1 & 0 & 0 & 1 &  \\ 2 & -1 & -1 & -1 & 1 \\
\end{array}\right)$$

The linear operator $\pi$ expresses the elementary row operations on the connectivity matrix $A$ needed to write it in lower triangular form; i.e. 
$$B^{t}A=\left(\begin{array}{ccccc} 2 &  &  &  &  \\ -1 & -1 &  &  &  \\ -1 & 0 & -1 &  &  \\ -1 & 0 & 0 & -1 &  \\ 1 & 1 & 1 & 1 & 1 \\
\end{array}\right)$$

Observe that the diagonal elements of $B^{t}A$ are the connectivity numbers $\alpha$; i.e. the coefficients of the $\overbrace{123}$ term in the $\pi$ relative to the basis $Part_{n}$ (the entries of the bottom row of the $B$ matrix).


\begin{prop}
The determinant of the connectivity matrix $A$ is $$det(A)=\prod_{\mathcal{A}\in Part_{n}}\alpha_{\mathcal{A}}$$
\end{prop}
\dem
Choose a coherent order in the basis $Part_{n}$ and consider the connectivity matrix $A$ relative to this order. Consider the matrix $B$ associated to the operator $\pi$ relative to the chosen coherent order in $Part_{n}$. We have that $B$ is a lower triangular matrix with ones in its diagonal, $$B=\left(
                         \begin{array}{cccc}
                           1 & 0 & \ldots & 0 \\
                           * & 1 & \ldots & 0 \\
                           \vdots & \vdots & \ddots & \vdots \\
                           * & * & \ldots & 1 \\
                         \end{array}
                       \right)$$ In particular, $$det(B)=1$$
As we argue before, we can think about the expression of the vector $\pi(\mathcal{A})$ in terms of the base $Part_{n}$ as elementary row operations on the matrix $A$ so, by Lemma \ref{lemaNueve}, we have the following identity: $$B^{t}A=\left(
                                                   \begin{array}{cccc}
                                                     \alpha_{\mathcal{O}_{1}}\ I_{\sharp\mathcal{O}_{1}} & 0 & \ldots & 0 \\
                                                     * & \alpha_{\mathcal{O}_{2}}\ I_{\sharp\mathcal{O}_{2}} & \dots & 0 \\
                                                     \vdots & \vdots & \ddots & \vdots \\
                                                     * & * & \ldots & \alpha_{\mathcal{O}_{k}}\ I_{\sharp\mathcal{O}_{k}} \\
                                                   \end{array}
                                                 \right)$$ Because $det(B^{t})=det(B)=1$ we have
                                                 $$det(A)=det(B^{t}A)= \prod_{i=1}^{k}\alpha_{\mathcal{O}_{i}}\ ^{\sharp\mathcal{O}_{i}}= \prod_{\mathcal{A}\in Part_{n}}\alpha_{\mathcal{A}}$$ and the proof is complete.
\fdem

\subsection{Step two: Connectivity numbers and reliability polynomial}

Although the connectivity matrix was motivated from a reliability problem, it is very interesting that the connectivity matrix determinant turns out to be a reliability calculation; i.e. The combinatorics of the connectivity numbers is encoded in the combinatorics of a reliability polynomial.

\begin{defi}
Denote by $K_{n}$ the graph with $n$ nodes and one edge joining every pair of nodes. Denote by $K_{n}^{\mathcal{A}}$ the resulting graph from the identification of the nodes $\{1,2,\ldots n\}$ in $K_{n}$ by the classes of the partition $\mathcal{A}$.
\end{defi}

\begin{lema}\label{CalculoNumerosConexidad}
Consider a partition $\mathcal{A}$. Then, $$R(K_{n}^{\mathcal{A}})= (-1)^{g} \alpha_{\mathcal{A}}\ p^{g} + \ldots$$ where $\alpha_{\mathcal{A}}$ is the connectivity number of the partition $\mathcal{A}$ and the dots denote terms of lower degree of the reliability polynomial.
\end{lema}
\dem
We claim that $$\pi(\mathcal{A})=\mathcal{A}\cdot\left(\prod_{\tau\ transp.\ /\ \langle\tau\rangle\nless\mathcal{A}}(e-\langle\tau\rangle)\right)$$
In effect, consider a connectivity state $\mathcal{B}$ such that $\mathcal{B}\nless \mathcal{A}$. There is a transposition $\tau$ in $\mathcal{B}$ not belonging to $\mathcal{A}$. Because $\mathcal{A}\subset\mathcal{A}\cdot\mathcal{B}$ and $\langle\tau\rangle\subset\mathcal{A}\cdot\mathcal{B}$ we have that $\mathcal{A}\cdot\langle\tau\rangle\subset\mathcal{A}\cdot\mathcal{B}$ so $$\mathcal{A}\cdot\mathcal{B}\cdot(\mathcal{A}-\mathcal{A}\cdot\langle\tau\rangle)=\mathcal{A}\cdot\mathcal{B}-\mathcal{A}\cdot\mathcal{B}=0$$ This implies that $$\mathcal{A}\cdot\mathcal{B}\cdot\left(\prod_{\tau\ transp.\ /\ \langle\tau\rangle\nless\mathcal{A}}(\mathcal{A}-\mathcal{A}\cdot\langle\tau\rangle)\right)=0$$ and we conclude that the factors $(\mathcal{A}-\mathcal{A} \cdot \mathcal{B})$ such that $\mathcal{B}\nless \mathcal{A}$ do not contribute in the original definition of $\pi$; i.e. Considering the transpositions not contained in $\mathcal{A}$ is enough. This way we have the following expression: $\pi(\mathcal{A})$:$$\pi(\mathcal{A})=\prod_{\tau\ transp.\ /\ \langle\tau\rangle\nless\mathcal{A}}(\mathcal{A}-\mathcal{A}\cdot\langle\tau\rangle)=\mathcal{A}\cdot\left(\prod_{\tau\ transp.\ /\ \langle\tau\rangle\nless\mathcal{A}}(e-\langle\tau\rangle)\right)$$ which proves the claim. In particular, the last expression implies that $$\alpha_{\mathcal{A}}=C_{0}-C_{1}+C_{2}-C_{3}+\ldots$$ where $C_{i}$ is the number of subsets $F$ with cardinality $i$ of the set of transpositions $\{\tau\ transp.\ /\ \langle\tau\rangle\nless\mathcal{A}\}$ such that $\langle F\rangle\cdot\mathcal{A}= S_{n}$.

Identifying the transposition $(ij)$ with the edge joining the nodes $i$ and $j$ of $K_{n}$, it is clear that $C_{i}$ is the pathsets number (operational states number) of $K_{n}^{\mathcal{A}}$ with just $i$ operational edges. This way we have that

\begin{eqnarray*}
  R(K_{n}^{\mathcal{A}}) &=& C_{0}(1-p)^{g} + C_{1}p(1-p)^{g-1} +C_{2}p^{2}(1-p)^{g-2} + \ldots \\
   &=& (-1)^{g}(C_{0}-C_{1}+C_{2}-C_{3}+\ldots)p^{g} +\ldots  \\
   &=& (-1)^{g}\alpha_{\mathcal{A}}\ p^{g}+\ldots
\end{eqnarray*} where the dots denote terms of lower degree.
\fdem

\begin{figure}
\begin{center}
  \includegraphics[width=0.7\textwidth]{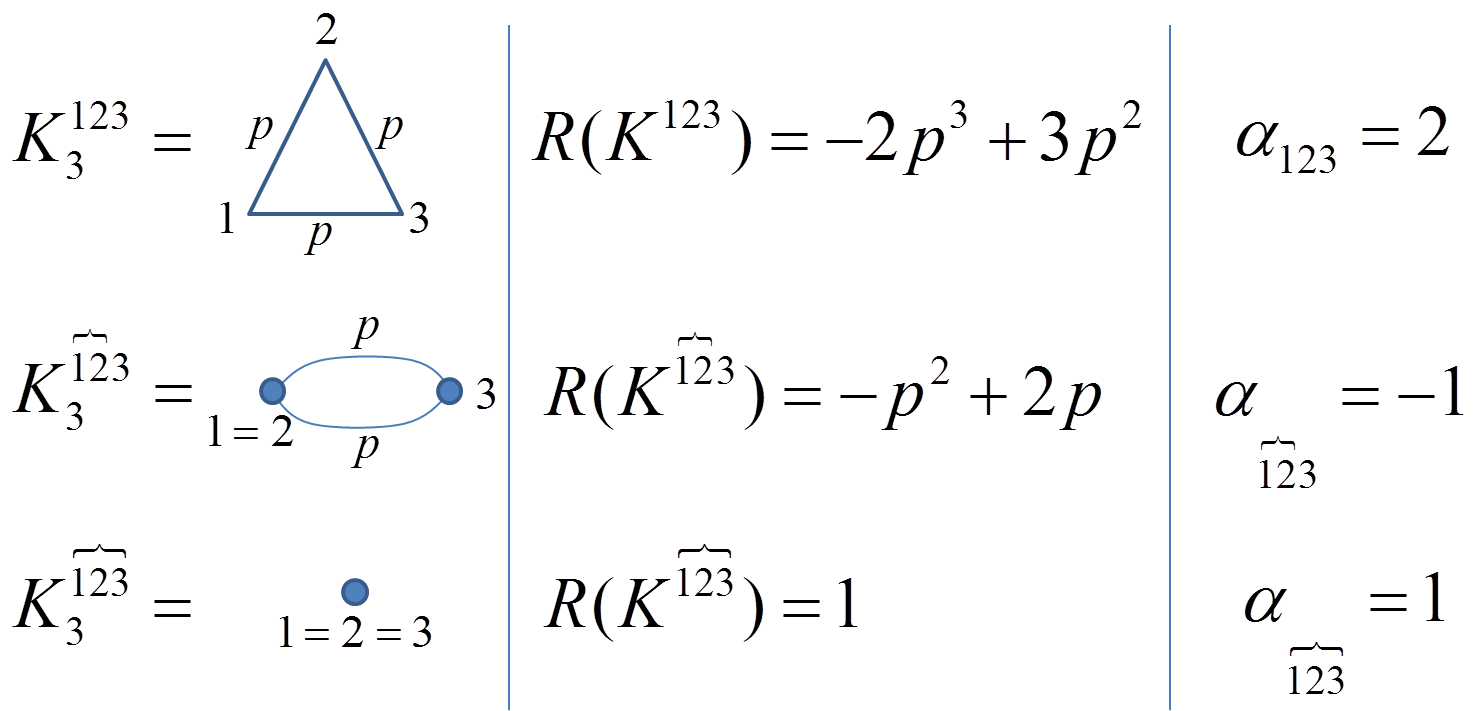}\\
  \end{center}
  \caption{Step 2 of the proof}\label{ExK3}
\end{figure}

Figure \ref{ExK3} shows an example in the $n=3$ case (irrelevant edges were not drawn).

\subsection{Step three: Connectivity numbers calculation}

This section is devoted to calculation of the highest degree term coefficient of the reliability polynomial of the graph $K_{n}^{\mathcal{A}}$. We assume the reader is aware of the simple factorization theorem ([Mo],[SC],[Co], [Sh]). Although it is not the main reference, it can also be found in the preliminaries of the previous first part paper [BR].

Let $G$ be a graph and consider its reliability polynomial $R(G)$. Denote by $mgr(R(G))$ the term of $R(G)$ whose degree equals the edge number of $G$. Observe that if $G$ has an irrelevant edge, then $mgr(R(G))=0$ and in case $mgr(R(G))$ is non zero, then this term equals the highest degree term of the polynomial. The following trick will be extremely useful in the following calculations.

\begin{lema}\label{Truco}
Consider a graph $G$ with $k$ edges between a pair of distinct nodes $i$ and $j$ of $G$. Consider the resulting graph $\tilde{G}$ by deleting $k-1$ edges between the nodes $i$ and $j$ of $G$. Then, $$mgr(R(G))=(-p)^{k-1}\ mgr\left(R(\tilde{G})\right)$$
\end{lema}
\dem
The result is clear for $k=1$. Suppose there are $k>1$ edges between the nodes $i$ and $j$ of $G$ and that the result holds for an amount less than or equal to $k-1$ of them. Consider an edge $a$ between the nodes $i$ and $j$. A simple factorization on the edge $a$ gives $$R(G)=p\ R(G\cdot a)+(1-p)R(G-a)$$ where $G\cdot a$ is the resulting graph by the contraction of $a$ and  $G-a$ is the resulting graph by deleting $a$. Observe that the edge number of $G\cdot a$ and $G-a$ is the edge number of $G$ minus one and because $k>1$, $G\cdot a$ has irrelevant edges. This way, $$mgr(R(G))=(-p)\ mgr(R(G-a))$$ By the inductive hypothesis, we get the result.
\fdem

The following is a well known corollary of Gilbert's formula [Gi],[Co] but we use our method to arrive at the same result.

\begin{lema}
The reliability polynomial of $K_{n}$ is $$R(K_{n})=\pm(n-1)!\ p^{g} + \ldots$$ where the highest degree of the expression $g$ equals the edge number of $K_{n}$: $$g=\left(
                                                                 \begin{array}{c}
                                                                   n \\
                                                                   2 \\
                                                                 \end{array}
                                                               \right)$$
\end{lema}
\dem
In this proof we will make an abuse of notation identifying the reliability polynomial with its graph. We claim that $$mgr(K_{n+1})=(-1)^{n+1}n\ mgr(K_{n})\ p^{n}$$ Because $K_{2}=p$ and $K_{1}=1$ we have the result for the $n=1$ case. Suppose the claim is true for every natural number less than or equal to $n$.

\begin{figure}
\begin{center}
  \includegraphics[width=0.35\textwidth]{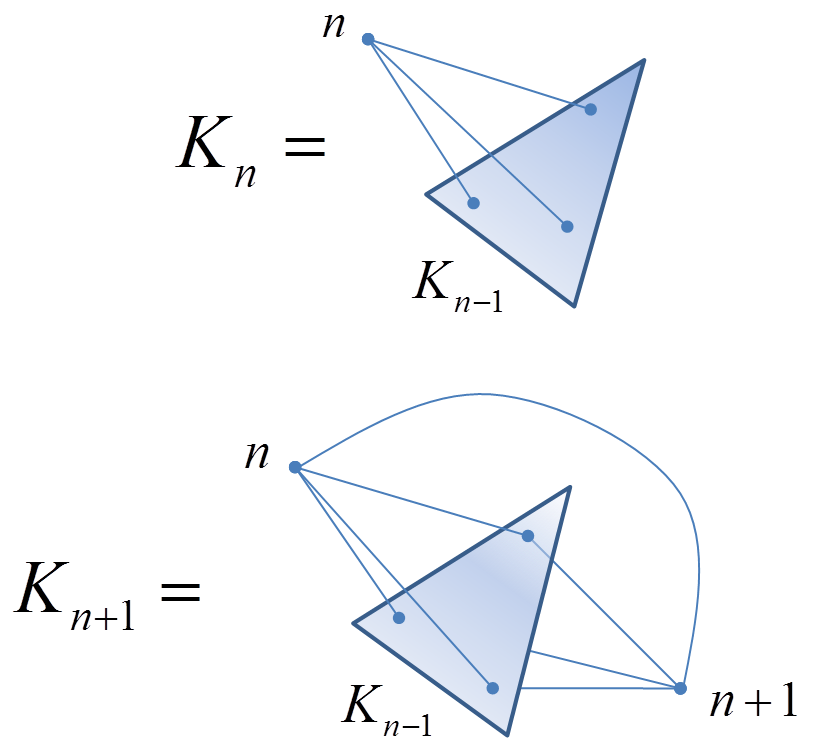}\\
  \end{center}
  \caption{Relation between the graphs $K_{n+1}$, $K_{n}$ and $K_{n-1}$}\label{RelacionGamas}
\end{figure}

\begin{figure}
\begin{center}
  \includegraphics[width=0.35\textwidth]{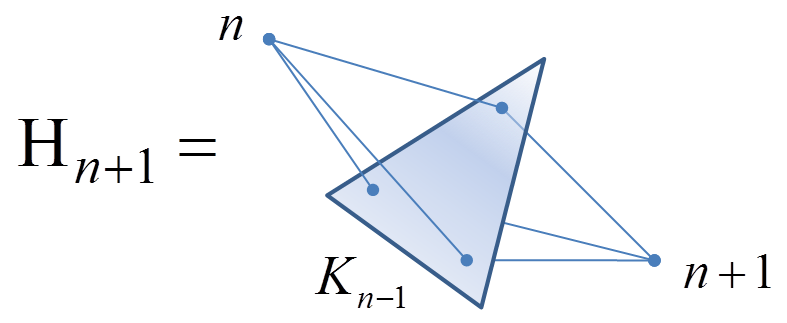}\\
  \end{center}
  \caption{The graph $H_{n+1}$}\label{GrafoH}
\end{figure}

The Figure \ref{RelacionGamas} shows the relation between the graphs $K_{n+1}$, $K_{n}$ and $K_{n-1}$. By a simple factorization on the edge joining the nodes $n$ and $n+1$ of the graph $K_{n+1}$ and the above lemma we have that $$K_{n+1}= p(-p)^{n-1}\ K_{n}+\ldots + (1-p)\ H_{n+1}$$ where the dots denote terms whose degree is less than the edge number of $K_{n+1}$ and the graph $H_{n+1}$ results from deleting the edge joining the nodes $n$ and $n+1$ of the graph $K_{n+1}$, see Figure \ref{GrafoH}.

By the inductive hypothesis, $$mgr(K_{n})=(-1)^{n}(n-1)\ mgr(K_{n-1})\ p^{n-1}$$ and the fact that the relation between the graphs $K_{n}$ and $K_{n-1}$ is the same as the one between $H_{n+1}$ and $K_{n}$ (both relations are the one point extension from the graph $K_{n-1}$, see Figures \ref{RelacionGamas} and \ref{GrafoH}), we have the following relation $$mgr(H_{n+1})=(-1)^{n}(n-1)\ mgr(K_{n})\ p^{n-1}$$ Then, we have that $$K_{n+1}= p(-p)^{n-1}\ K_{n}+\ldots + (1-p)(-1)^{n}(n-1)\ K_{n}\ p^{n-1}+ \ldots= (-1)^{n+1}n\ K_{n}\ p^{n}+\ldots$$ where the dots denote terms whose degree is less than the edge number of $K_{n+1}$. We conclude that $$mgr(K_{n+1})=(-1)^{n+1}n\ mgr(K_{n})\ p^{n}$$ which proves the claim. This recursive relation shows that $mgr(K_{n})$ is non zero so it equals the highest degree term of the reliability polynomial for $K_{n}$: $$K_{n}=(-1)^{n+(n-1)+\ldots 2}(n-1)!p^{(n-1)+\ldots 1}+ \ldots$$ and this concludes the lemma.
\fdem

\begin{lema}\label{CalculoMagico}
Consider a partition $\mathcal{A}$ with $m$ classes: $\mathcal{A}=\{ a_{1}, a_{2}, \ldots a_{m}\}$. Then, the reliability polynomial of $K_{n}^{\mathcal{A}}$ is:
$$R(K_{n}^{\mathcal{A}})= \pm (m-1)!\ p^{g} + \ldots$$ where the highest degree $g$ of the expression equals edge number of $K_{n}^{\mathcal{A}}$ after taken out the irrelevant edges: $$g=\sum_{i\neq j}(\sharp a_{i})\ (\sharp a_{j})$$
\end{lema}
\dem
We will make the same abuse we did in the proof before identifying the reliability polynomial with its graph and we also will identify the class $a_{i}$ with its cardinality $\sharp a_{i}$ for notational convenience. The graph $K_{n}^{\mathcal{A}}$ has $m$ nodes (these are the $m$ classes of $\mathcal{A}$), $\left(
                        \begin{array}{c}
                          a_{i} \\
                          2 \\
                        \end{array}
                      \right)$ irrelevant edges in each node $i$ respectively and $a_{i}a_{j}$ edges joining the nodes $i$ and $j$.

Consider the graph $\bar{K}_{n}^{\mathcal{A}}$ resulting from deleting all the irrelevant edges of the graph $K_{n}^{\mathcal{A}}$. This way $K_{n}^{\mathcal{A}}$ and $\bar{K}_{n}^{\mathcal{A}}$ have the same reliability polynomial. By the Lemma \ref{Truco} we have the following relation between the graphs $\bar{K}_{n}^{\mathcal{A}}$ and $K_{m}$:

\begin{eqnarray*}
  mgr(\bar{K}_{n}^{\mathcal{A}}) &=& (-p)^{\sum_{i\neq j}(a_{i}a_{j}-1)}\ mgr(K_{m}) \\
   &=& (-p)^{\sum_{i\neq j}(a_{i}a_{j}-1)}\ (\pm (m-1)!)\ p^{\left(
                                                             \begin{array}{c}
                                                               m \\
                                                               2 \\
                                                             \end{array}
                                                           \right)} \\
   &=&  \pm (m-1)!\ p^{\sum_{i\neq j}a_{i}a_{j}}
\end{eqnarray*} and this concludes the proof.
\fdem

This concludes the third step of the proof. We have finally proved the determinant formula:

\begin{prop}\label{AInvertible}
The determinant of the connectivity matrix $A$ is: $$det(A)= \pm \prod_{\mathcal{A}\in Part_{n}}(m_{\mathcal{A}}-1)!$$ where $m_{\mathcal{A}}$ is the number of classes in the partition $\mathcal{A}$. In particular, $A$ is invertible.
\end{prop}

As an example, consider the $n=4$ case. Recall the partitions conjugation classes are:

\begin{eqnarray*}
  \mathcal{O}_{1} &=& \{\ 1234\ \} \\
  \mathcal{O}_{2} &=& \{\ 12 \overbrace{34}, 13\overbrace{24}, 23\overbrace{14}, 1\overbrace{23}4, \overbrace{13}24, \overbrace{12}34\ \} \\
  \mathcal{O}_{3} &=& \{\ \overbrace{14}\overbrace{23}, \overbrace{13}\overbrace{24}, \overbrace{12}\overbrace{34}\ \} \\
  \mathcal{O}_{4} &=& \{\ 1\overbrace{234}, 2\overbrace{134}, 3\overbrace{124},\overbrace{123}4\ \} \\
  \mathcal{O}_{5} &=& \{\ \overbrace{1234}\ \}
\end{eqnarray*} then we have that $m_{\mathcal{O}_{i}}$ equals $4, 3, 2, 2, 1$ and the partitions conjugation classes cardinalities are $1, 6, 3, 4, 1$ respectively. By the above proposition the determinant of the connectivity matrix is: $$det(A)=\pm(4-1)!\ ^{1}(3-1)!\ ^{6}(2-1)!\ ^{3}(2-1)!\ ^{4}(1-1)!\ ^{1}= \pm 384$$

\section{Acknowledgements}
The author is grateful to Franco Robledo and the anonymous referees for their careful reading and valuable suggestions in the improvement of the paper.


\end{document}